\documentclass[dvips,aap]{imsart}

\RequirePackage[OT1]{fontenc}
\RequirePackage{amsthm,amsmath,natbib,amssymb,amsfonts,graphicx}
\RequirePackage[colorlinks]{hyperref}
\RequirePackage{hypernat}
\RequirePackage{comment}

\startlocaldefs
\numberwithin{equation}{section}
\theoremstyle{plain}
\newtheorem{thm}{Theorem}[section]
\newtheorem{lem}[thm]{Lemma}
\newtheorem{cor}[thm]{Corollary}
\newtheorem{prop}[thm]{Proposition}
\theoremstyle{definition}
\newtheorem{defi}[thm]{Definition}

\theoremstyle{remark}

\newtheorem{rem}[thm]{Remark}

\newcommand{\E}{\mathbb{E}}

\endlocaldefs

\begin{document}

\begin{frontmatter}
\title{Holomorphic transforms with application to affine processes}
\runtitle{Holomorphic transforms}

\begin{aug}
\author{\fnms{Denis} \snm{Belomestny}\thanksref{t1}\ead[label=e1]{belomest@wias-berlin.de}},
\author{\fnms{J\"org} \snm{Kampen}\thanksref{t2}\ead[label=e2]{kampen@wias-berlin.de}},
\and
\author{\fnms{John} \snm{Schoenmakers}\thanksref{t2}
\ead[label=e3]{schoenma@wias-berlin.de}
\ead[label=u1,url]{http://www.wias-berlin.de/~schoenma}}

\thankstext{t1}{Partially supported by  SFB 649 `Economic Risk'.}
\thankstext{t2}{Partially supported by  DFG Research Center \textsc{Matheon} `Mathematics for Key
Technologies' in Berlin}
\runauthor{D. Belomestny et al.}

\affiliation{Weierstrass Institute}

\address{Weierstrass Institute for Applied Analysis and Stochastics\\
Mohrenstr. 39
\\
10117, Berlin, Germany
\\
\printead{e1}
\\
\phantom{E-mail:\ }\printead*{e2}
\\
\phantom{E-mail:\ }\printead*{e3}
}

\end{aug}

\begin{abstract}
In a rather general setting of It\^o-L\'evy processes we study a class of
transforms (Fourier for example)  of the state variable of a process which are holomorphic
in some disc around time zero in the complex plane.
 We show that such  transforms are related to a system of
analytic vectors for the generator of the process, and we state
conditions which allow for  holomorphic extension of these transforms into a strip
which contains the positive real axis. Based on these extensions we develop
a functional series expansion of these transforms in terms of the constituents of the generator.
As application, we  show that
for multidimensional affine It\^o-L\'evy processes with state
dependent jump part
the Fourier transform is
holomorphic in a time
strip under some stationarity conditions, and give
log-affine series representations for the transform.
\end{abstract}

\begin{keyword}[class=AMS]
\kwd[Primary ]{60J25}
\kwd[; secondary ]{91B28}
\end{keyword}

\begin{keyword}
\kwd{It\^o-L\'evy processes}
\kwd{holomorphic transforms}
\kwd{affine processes}
\end{keyword}

\end{frontmatter}

\section{Introduction}
Transforms are an important tool in the theory of (ordinary and partial)
differential equations and in  stochastic analysis.
In probability theory the Fourier transform  of a random variable, which
represents  the characteristic function of the
corresponding distribution, is widely used. Fourier (and Laplace)
transforms have
become increasingly popular in
mathematical  finance as well.  On the one hand, for example, via a complex
Laplace transform and the convolution theorem one may derive pricing formulas
for one dimensional European options (e.g., see \cite{EP}). On the other hand,
Laplace and
Fourier transforms   are known in  closed form for many classes of processes.
A famous
example is the so-called L\'evy-Khintchine formula which provides an explicit
expression for the characteristic function of a L\'evy process. More recent
financial
literature goes well beyond L\'evy processes and attempts to establish explicit or
semi-explicit formulas for derivatives where underlyings are modelled,
for instance, by affine processes
(see among others \cite{DPS},\cite{F} and \cite{DFS}) or
more general It\^o-L\'evy processes (e.g. see \cite{OS}). Theoretical analysis of affine
processes is done in the seminal paper by Duffie, Filipovi\'c and Schachermayer
\cite{DFS} and has led to a unique characterization of affine processes. In particular, it
is shown that
the problem of determining  the (conditional) Fourier transform of  an
affine process $X_s$ corresponds to the problem of solving a system of
generalized Riccati differential equations  in the time variable $s$ (see \cite{DFS}).
Although closed form solutions of this system can be found in important cases,
there is no generic approach to solve such a system in the general
multi-dimensional case.
In this article we establish some kind of functional series
representation for the Fourier transform, hence the characteristic function of the process
under consideration and, in principle, for more general transforms.
The most  natural one is  a Taylor expansion in time $s$ around $s_0=0$. This approach leads to
expansions where the coefficients can be recursively computed without solving linear or
non-linear differential equations as, for instance, in the case of the so called WKB expansions
for transition densities   (e.g. see \cite{K}).
Unfortunately, it turns out that in many cases the resulting
power series converges only in $s=s_0=0.$
This problem corresponds to a difficulty which is well known in
semi-group theory and in the theory of parabolic differential equations:
small time expansions for the solutions of parabolic
equations are usually possible in a neighborhood of some  $s_0 >0,$
while an expansion  around $s_0=0$ may have zero convergence radius.
In this paper we prove that for a generator with affine coefficients
the Fourier transform extends holomorphically into a disc around $s_0=0$
and a strip containing the positive real axes,  under some mild regularity conditions.
Then, for multi-dimensional affine processes we obtain
convergent expansions for the Fourier transform  and its logarithm on the whole time line.
Hence, we  have (affine) series representations for the exponent of the characteristic
function  of a  general multi-dimensional affine process.
More generally, we develop a framework
based  on a concept of analytic vectors   which allows for functional series
expansions for a class of  holomorphic transforms which covers the standard Fourier transform.

The outline of the paper is as follows. The basic setup is described in
Section~\ref{BS}. In Section \ref{AVT} we introduce the notion of analytic vectors
associated with a
given generator and study functional series expansions for the corresponding
transform.
Section \ref{genA} is devoted
to the Cauchy problem for affine generators.
In Section \ref{LAFF} we derive series representations for the logarithm
of the Fourier transform  corresponding to a generator with affine coefficients.
Section~\ref{EX} gives an explicit representation for these expansions in a
one-dimensional case.
Finally, Section \ref{ItoLevy} contains results for affine
It\^o-L\'evy processes which mainly follow from  previous sections.
More technical proofs are given
in the Appendix.


\section{Basic setup}
\label{BS}
Let $(\Omega, \mathcal{F}, (\mathcal{F}_{t})_{t\ge0},P)$ be a standard filtered probability space
where the filtration $(\mathcal{F}_{t})$
satisfies `the usual conditions'.  On this space we consider for each $x\in\mathbb{R}^{n}$ a compensated Poisson random
measure $\widetilde{N}(x,dt,dz,\omega)=N(x,dt,dz,\omega)-v(x,dz)dt$ on
$\mathbb{R}_{+}\times\mathbb{R}^{n},$ where $N$ is a Poisson random measure
with (deterministic) intensity kernel of the form $v(x,dz)dt=\mathbb{E}\, N(x,dt,dz)$ satisfying
$v(x,B)<\infty$ for any $B\in\mathcal{B}(\mathbb{R}^{n})$ such that
$0\notin\overline{B}$ (closure of $B$). Hence $N$ is determined by
\[
P\left[  N(x,(0,t],B)=k\right]  =\exp(-tv(x,B))\frac{t^{k}v^{k}(x,B)}%
{k!},\quad k=0,1,2,...
\]
In particular, for any $B\in\mathcal{B}(\mathbb{R}^{})$ with $0\notin$ $\overline{B}$ and $x\in\mathbb{R}^n,$ the process $M_{t}^{B,x}:=\widetilde{N}(x,(0,t],B)$ is a (true) martingale.
Further, by assumption, the kernel $v$ satisfies,%
\[
v(x,\{0\})=0,\quad\int_{\mathbb{R}^{n}}(|z|^{2}\wedge
|z|)v(x,dz)<\infty,\quad x\in\mathbb{R}^{n}.
\]

Let us assume $W(t)$ to be a standard Brownian motion  in $\mathbb{R}^{m}$ living on our basic probability space, and consider
the It\^o-L\'evy SDE:
\begin{equation}
dX_{t}=b(X_{t})dt+\sigma(X_{t})dW(t)+\int_{\mathbb{R}^{n}} z\widetilde{N}(X_{t-},dt,dz),\quad X_0=x, \label{ItoLev}%
\end{equation}
for  deterministic functions $b:\mathbb{R}%
^{n}\rightarrow\mathbb{R}^{n},$ $\sigma:\mathbb{R}^{n}\rightarrow
\mathbb{R}^{n}\times\mathbb{R}^{m},$  which satisfy sufficient regularity and/or  mutual consistency conditions
such that  (\ref{ItoLev}) has a unique strong solution
$t\rightarrow X_t,$  called an It\^o-L\'evy process,
which can be regarded as a strong Markov process  (e.g., see \cite{P},  \cite{OS}).

As a well-known fact, the above process $X$ can be  connected  to some kind of evolution equation in a  natural way.
%
%
In this context we consider a 'pseudo generator'
\begin{equation}
A^\sharp:\mathcal{D}(A^\sharp)\subset \mathcal{C}^{(2)}\subset \mathcal{C}
\longrightarrow\mathcal{C}, \label{PSG}
\end{equation}
where $\mathcal{C}:=C(\mathbb{R}^{n})$ is the space of  continuous functions $f:$ $\mathbb{R}^{n}\rightarrow \mathbb{C},$ equipped with the topology of uniform convergence on compacta, and $\mathcal{C}^{(2)}$ is the space of functions $f$ $\in$ $\mathcal{C}$ which are two times continuously differentiable. Further, $f\in\mathcal{D}(A^\sharp)$ iff $f\in\mathcal{C}^{(2)},$ and
\begin{align}
A^\sharp f(x)  &  :=\frac{1}{2}\sum_{i,j=1}^{n}a_{ij}(x)\frac{\partial^{2}f}{\partial
x_{i}\partial x_{j}}+\sum_{i=1}^{n}b_{i}(x)\frac{\partial f}{\partial x_{i}%
}\label{Gen}\\
&  +\int_{\mathbb{R}^{n}}\left[  f(x+z)-f(x)-\frac{\partial
f}{\partial x}\cdot z\right]  v(x,dz),\ \rm{with}\  a:=\sigma\sigma^{\top},\nonumber
\end{align}
exists and is such that
$A^\sharp f\in\mathcal{C}.$  In this respect we assume that the building blocks  $b(x),$ $\sigma(x),$ $v(x,dz)$ of SDE (\ref{ItoLev})
have bounded derivatives of any order with respect to $x.$ Clearly, $\mathcal{D}(A^\sharp)$ is dense in $\mathcal{C}$ and it can be shown that the operator $A^\sharp$ thus defined is closable. The closure of $A^\sharp$  is denoted by
\begin{equation}\label{clA}
A:\mathcal{D}(A)\subset \mathcal{C}
\longrightarrow\mathcal{C}.
\end{equation}
As such $A$ can be seen as a relaxation of the notion of a generator of a strongly continuous Feller-Dynkin semigroup
associated with the process $X,$ for which the Hille-Yosida theorem applies. This semigroup  is usually defined on the Banach space $C_0(\mathbb{R}^{n}),$ i.e. the set of  continuous functions on $\mathbb{R}^{n}$ which vanish at infinity, equipped with the supremum norm, and its generator coincides literally with $A^\sharp$ on a dense subdomain of $C_0(\mathbb{R}^{n}).$ By slight abuse of terminology however, we will also refer to $A$
as 'generator' when $A$ is considered in connection with the process $X$ given via (\ref{ItoLev}).
Let now $\mathfrak{F}:=\{f_{u},$ $u\in I\}$ $\subset$ $\mathcal{C}$ be a dense subset of bounded continuous  functions
$f_{u}:\mathbb{R}^{n}\rightarrow\mathbb{C}$ which have bounded derivatives of any order.
With respect to the (closed) generator  (\ref{clA}) we  consider  for each $f_u$ $\in$ $\mathfrak{F}$ the  (generalized) Cauchy problem
\begin{eqnarray}\label{FBK1}
 \frac{\partial\widehat{p}}{\partial s}(s,x,u)&=& A\widehat{p}\, (s,x,u),
\\
\nonumber
\widehat{p}(0,x,u)&=&f_{u}(x),\quad s\ge0,\quad x\in \mathfrak{X}\subset\mathbb{R}^{n},
\end{eqnarray}
where $\mathfrak{X}$ is some open (maximal) domain, and assume that problem (\ref{FBK1}) has a unique solution. For instance, see \cite{AM} for  mild conditions which guarantee existence
of unique global solutions of integro-differential evolution problems.
In particular, if some ellipticity condition is satisfied we may have $\mathfrak{X}$ $=$ $\mathbb{R}^{n}.$
For mixed type generators, such as  affine generators,
existence, uniqueness, and the maximal domain $\mathfrak{X}$ has to be considered case by case. For example, if
$b(x),$ $(\sigma\sigma^{\top})(x),$ $v(x,dz)$  are affine in $x$ and satisfy the admissibility conditions in \cite{DFS}, existence and uniqueness are ensured in a  domain  of the form $\mathfrak{X}$ $=$ $\mathbb{R}^{l}\times\mathbb{R}^{n-l}_+\subset\mathbb{R}^{n}$ (for details see \cite{DFS}). We underline, however, that in this article the main focus is on  functional series representations for the solution of \eqref{FBK1}, and  we therefore  merely assume that sufficient regularity conditions
for the coefficients in  (\ref{Gen}) (hence (\ref{ItoLev}))
are fulfilled.
\begin{rem}\label{AX}{\rm
In our analysis we often consider the pseudo generator (\ref{Gen}) and it's closure (\ref{clA}) on $\mathcal{C}:=C(\mathfrak{X})$, for an open domain \(\mathfrak{X}\subset \mathbb{R}^{n}\), rather than $C(\mathbb{R}^{n}).$
For notational convenience (while slightly abusing notation) we will denote these respective
operators with $A^\sharp$ and $A$ also.
}
\end{rem}

  If $A$ is the generator of the process (\ref{ItoLev}), the solution $\widehat{p}(s,x,u)$  has the probabilistic representation%
\[
\widehat{p}(s,x,u)=\mathbb{E}\left[\,f_{u}(X_{s}^{0,x})\right],
\]
where $X^{0,x}$ is the unique strong solution of (\ref{ItoLev}) with $X_{0}^{0,x}=x$
a.s. We refer to $\widehat{p}(s,x,u)$ as  \textit{generalized transform} of
the process $X_{s}^{0,x}$ associated with $\mathfrak{F}$. As a canonical
example we may consider
\begin{equation}
f_{u}(x):=e^{\mathfrak{i}u^{\top}x},\quad u\in\mathbb{R}^{n},\label{SF}%
\end{equation}
in which case (\ref{FBK1}) yields the characteristic function $\widehat
{p}(s,x,u)=\E[e^{\mathfrak{i}u^{\top}X_{s}^{0,x}}].$

By using multi-index notation, the integral term in (\ref{FBK1}) may be formally expanded as%
\begin{align*}
&  \int_{\mathbb{R}^{n}}\left[  f(x+z)-f(x)-\frac{\partial
f}{\partial x}\cdot z\right]  v(x,dz)=\\
&
 \sum_{|\alpha|\geq2}\frac{1}{\alpha!}\partial_{x^{\alpha}}f(x)\int
z^\alpha v(x,dz)
  =:\sum_{|\alpha|\geq2}\frac{1}{\alpha!}m_{\alpha}(x)\partial_{x^{\alpha}%
}f.
\end{align*}
Hence, we may write formally the generator as an infinite order differential operator
\begin{equation} A=\sum_{|\alpha|>0}\mathfrak{a}_{\alpha}(x)\partial_{x^{\alpha}}\label{FBK2}
\end{equation}
with obvious definitions of the coefficients $\mathfrak{a}_{\alpha}(x)$ for $\alpha>0.$

\section{Analytic vectors and transforms}\label{AvT}
\label{AVT}
First we introduce
the notion of a set of analytic vectors associated with  an operator $A.$

\begin{defi}\label{AV} $\mathfrak{F}=\{f_{u},$ $u\in I\}$ is a set of analytic vectors for an operator $A$ in an open region $\mathfrak{X,}$  if
\begin{description}
\item[\rm{(i)}] $A^{k}f_{u}$ exists for any $u\in I$  and  $k\in \mathbb{N},$

\item[\rm{(ii)}] for every $u\in I$ there exists  $R_{u}>0$
such that for all $x\in\mathfrak{X,}$
\begin{equation*}\label{exp}
\sum_{k=0}^{\infty}\frac{s^{k}}{k!}\left\vert A^{k}f_{u}(x)\right\vert
<\infty,\quad0\leq s<R_{u},
\end{equation*}
where the convergence 
is uniform in $x$ over any compact subset of
$\mathfrak{X.}$
\end{description}
\end{defi}

Thus, if $\mathfrak{F}$\ is a set of analytic vectors in the sense of Definition~\ref{AV} then
for all $x\in\mathfrak{X}$ the map
\begin{equation} \label{ConvP}
s\rightarrow P_{s}f_{u}(x):=\sum_{k=0}^{\infty}\frac{s^{k}}{k!}A^{k}%
f_{u}(x), \quad \left\vert s\right\vert <R_{u}
\end{equation}
is holomorphic in the complex disc $D_{0}:=\left\{  s\in\mathbb{C}:\left\vert
s\right\vert <R_{u}\right\}  $ and the series converges uniformly in $x$
over any compact subset of $\mathfrak{X.}$
In fact, \( P_{s}f_{u}(x) \) coincides with \( \widehat{p}(s,x,u) \) for
\(0\le s<R_{u}\).

\begin{prop}\label{propa}$\ $ If $\mathfrak{F}$\ is a set of analytic vectors, the map $(s,x)\rightarrow P_{s}f_{u}(x)$
defined in {\rm (\ref{ConvP})} satisfies {\rm (\ref{FBK1})} for all $s,$ $\left\vert
s\right\vert <R_{u}$ and $x\in\mathfrak{X.}$ In particular we have $P_{s}f_{u}(x)=\widehat{p}(s,x,u)$, $0\le s<R_{u}.$
\end{prop}

\begin{proof} Obviously, $P_{0}f_{u}(x)=f_{u}(x)$ for $x\in\mathfrak{X.}$
Set (see Remark~\ref{AX})
\[
P_{s}^{(N)}f_{u}(x):=\sum_{k=0}^{N}\frac{s^{k}}{k!}A^{k}f_{u}(x),
\]
 then both $P_{s}^{(N)}f_{u}(x)$ and
\[
AP_{s}^{(N)}f_{u}(x):=\sum_{k=0}^{N}\frac{s^{k}}{k!}A^{k+1}f_{u}(x)
\]
converge uniformly for any $x$ in a compact subset of $\mathfrak{X}$ and for any $s$
satisfying $\left\vert s\right\vert <R_{u}-\varepsilon$ with  arbitrary small \( \varepsilon \). Hence,
since $A$ is closed,%
\[
AP_{s}f_{u}(x)=\sum_{k=0}^{\infty}\frac{s^{k}}{k!}A^{k+1}f_{u}(x)=\frac
{\partial}{\partial s}\sum_{k=0}^{\infty}\frac{s^{k}}{k!}A^{k}f_{u}%
(x)=\frac{\partial}{\partial s}P_{s}f_{u}(x).
\]
\end{proof}

In order to study generalized transforms associated with a set of analytical vectors
$\mathfrak{F}$ in domains containing the non-negative real axis we introduce for $\eta>0$ the sequence
\begin{align}
q_{k}^{(\eta)}(x,u) &  :=\frac{1}{k!}\prod\limits_{r=0}^{k-1}
(\eta^{-1}A+rI)f_{u}(x)
\label{qprod}\\
&  =:\frac{1}{k!}\sum_{r=0}^{k}c_{k,r}\eta^{-r}A^{r}f_{u}(x),\quad
x\in\mathfrak{X,}\quad u\in I,\quad k=0,1,2,...\label{qseq}%
\end{align}
In (\ref{qseq}) the coefficients $c_{k,r},$ $0\leq r\leq k,$ are determined by the
identity
\[{\displaystyle\prod\limits_{r=0}^{k-1}}
(z+r)=z(z+1)\cdot\ldots\cdot(z+k-1)\equiv\sum_{r=0}^{k}c_{k,r}z^{r},
\]
and are usually called unsigned Stirling numbers of the first kind.
These numbers satisfy  $c_{0,0}=1$ and
\begin{align}
\nonumber
c_{k,0} &  \equiv0,\quad c_{k,k}\equiv1,\\
\label{ck}
c_{k+1,j} &  =kc_{k,j}+c_{k,j-1},\quad1\leq j\leq k,
\end{align}
if  $k\geq1$.
Obviously, the following recursion is equivalent to (\ref{qprod}),
\begin{equation}
(k+1)q^{(\eta)}_{k+1}(x,u)=\eta^{-1}Aq^{(\eta)}_{k}(x,u)+kq^{(\eta)}_{k}(x,u),\quad k\geq0,\quad\text{
}x\in\mathfrak{X}. \label{rec}%
\end{equation}
The next theorem provides a functional series representation for the solution of \eqref{FBK1} for all $s\ge0,$ under certain conditions.
\begin{thm}\label{analth}
Let $\mathfrak{F}$\ be a set of analytic vectors in
the sense of Definition~\ref{AV}, $u\in I$ be fixed, and the sequence $(q_{k}^{(\eta)})$ be
defined as in (\ref{qseq}). Let   $\widehat{p}$ be the solution of the Cauchy problem (\ref{FBK1}).
Then the following statements are equivalent:

$(i)$ There exists a constant $R_{u}>0$ such that for each $x\in\mathfrak{X},$
the map
$
s\rightarrow\widehat{p}(s,x,u)
$
 has a holomorphic extension to the domain%
\[
G_{R_{u}}:=\left\{  z:\left\vert z\right\vert <R_{u}\right\}
\cup\left\{  z:\operatorname{Re}z>0\ \wedge\ |\operatorname{Im}z|<R_{u}%
\right\},
\]
see Figure~\ref{Perform}.

\begin{figure}[pth]
\centering \includegraphics{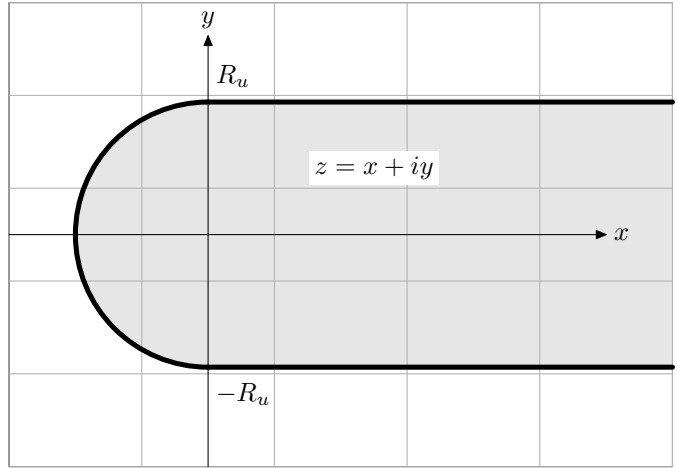}\caption{ Domain $G_{R_{u}} $ on
the complex plane }%
\label{Perform}%
\end{figure}

$(ii)$ There exists an $\eta_{u}>0$ such that  for each $x\in\mathfrak{X}$ the following
 series representation holds:
\[
\widehat{p}(s,x,u)=\sum_{k=0}^{\infty}q_{k}^{(\eta_{u})}(x,u)\left(
1-e^{-\eta_{u}s}\right)  ^{k},\quad 0\le s<\infty.%
\]
Moreover, the series converges uniformly for $(x,s)$ running through any compact subset of $\mathfrak{X} $ $\times$ $\left\{  s\in\mathbb{R}:0\le s< R_u\right\}.$

$(iii)$
The solution  $\widehat{p}$ of the Cauchy problem (\ref{FBK1}) is holomorphically extendable to
$[0,\infty),$
there exists  $\eta_{u}>0$  such that
\begin{equation}\label{convratio}
\overline{\lim}_{k\rightarrow\infty}\sqrt[k]{\left\vert q_{k}^{(\eta_{u}%
)}(x,u)\right\vert }\le1,\quad x\in \mathfrak{X},
\end{equation}
and, there exists  $\varepsilon_u,$ $0$ $<$ $\varepsilon_u$ $<$ $1,$ such that the
series
\begin{equation}\label{eps3}
\sum_{k=0}^{\infty}q_{k}^{(\eta_{u})}(x,u)w^{k}%
\end{equation}
converges uniformly for $(x,w)$ running through any compact subset of $\mathfrak{X} $ $\times$ $\left\{  w\in\mathbb{C}:\left\vert w\right\vert
<1-\varepsilon_u\right\}.$

\end{thm}

\begin{proof}
See Appendix.
\end{proof}

\begin{rem} The implication
 $(iii)'$ $\Rightarrow$ $(i),$ where statement $(iii)'$ consists of (\ref{convratio}), and (\ref{eps3}) with
$\varepsilon_u=0$ holds as well. That is, loosely speaking, if in $(iii)$  series (\ref{eps3}) converges uniformly on
all compact subsets of $\mathfrak{X} $ $\times$ $\left\{  w\in\mathbb{C}:\left\vert w\right\vert
<1\right\},$ the holomorphy assumption on $\widehat p$  can be dropped.
\end{rem}
\begin{rem}
In order to use the representation in \( (ii) \) one has to choose \( \eta_{u} \). In
fact, \( \eta_{u} \) can be related to \( R_{u} \) via \( \eta_{u}=\pi/R_{u} \) and hence
increases with decreasing  \( R_{u} \).
\end{rem}
It is important to note Theorem~\ref{analth} concerns the solution of the Cauchy problem
(\ref{FBK1}) connected with a general operator $A.$ In particular, all
criteria in this theorem are of pure analytic nature and via (\ref{qprod}), respectively
(\ref{qseq}), exclusively formulated in terms of the $A^k f_u (x)$, i.e.  coefficients
in Definition~\ref{AV}. In the case where
$A$ is the generator of a Feller Dynkin process
one can formulate a sufficient {\em probabilistic} criterion for Theorem~\ref{analth}{\em-(i)}:

\begin{prop}\label{probcrit}
Let \( \mathfrak{F} \) be a set of analytic
vectors  in the sense of Definition~\ref{AV} and let the Markov process \( \{X_{t}\} \) be  associated with the generator \( A \). If in addition,
for every $u\in I$ there
exists a radius $R_{u}$ such that for any \( t\geq 0 \)
\begin{equation}
\sum_{k=0}^{\infty}\frac{s^{k}}{k!}\left\vert A^{k} \E [f_{u}(X^{0,x}_{t})]\right\vert <\infty,\quad0\leq s<R_{u},\label{AV1}%
\end{equation}
uniformly in $x$ over any compact subset of $\mathfrak{X}$,
then  Theorem~\ref{analth}{\em-(i)} holds.
\end{prop}
The statement  is a direct consequence of the following ``quasi'' semi-group property of the transition
operator \( P_{t} \).
\begin{prop}
Let $\mathfrak{F}$\ be a set of analytic vectors satisfying (\ref{AV1}).
Then, for all $x\in\mathfrak{X}$ and
all $t\geq 0,$ the generalized transform \( \widehat p({t+s},x,u) \) can be represented as
\begin{equation}
 \widehat p({t+s},x,u)=\sum_{k=0}^{\infty}\frac{s^{k}}%
{k!}A^{k}\E[f_{u}(X_{t}^{0,x})]  ,\quad 0\leq s<R_{u}, \label{ConvP1}%
\end{equation}
where the series converges uniformly in $x$ over any compact subset of $\mathfrak{X}.$
\end{prop}
\begin{proof}
Denote the right-hand-side of \eqref{ConvP1} by \( {\widetilde p}(t,s,x,u) \).
Obviously, ${\widetilde p}(t,0,x,u)=\E[f_{u}(X_{t}^{0,x})].$ Set
\[
{\widetilde p}^{(N)}(t,s,x,u):=\sum_{k=0}^{N}\frac{s^{k}}{k!}A^{k}\E[f_{u}(X_{t}^{0,x})],
\]
then both ${\widetilde p}^{(N)}(t,s,x,u)$ and
\[
A{\widetilde p}^{(N)}(t,s,x,u)=\sum_{k=0}^{N}\frac{s^{k}}{k!}A^{k+1}\E[f_{u}(X_{t}^{0,x})]
\]
converge for $N\rightarrow\infty$ uniformly over any compact subset of
$\mathfrak{X,}$ and $s$ with $\left\vert s\right\vert $ $<R_{u}-\varepsilon,$ for
an arbitrary small $\varepsilon>0$. Hence, for $\left\vert s\right\vert $
$<R_{u}-\varepsilon,$ we have
\begin{align*}
\frac{\partial}{\partial
s}{\widetilde p}^{(N)}(t,s,x,u)&=\sum_{k=0}^{N-1}\frac{s^{k}}%
{k!}A^{k+1}\E[f_{u}(X_{t}^{0,x})]=A{\widetilde p}^{(N-1)}(t,s,x,u),\\
{\widetilde p}(t,0,x,u)  & =\widehat{p}(t,x,u)
\end{align*}
and thus, by closeness of the operator $A$ and uniqueness of the
Cauchy problem (\ref{FBK1})-(\ref{FBK2}), we have \( {\widetilde p}(t,s,x,u)=\widehat{p}(t+s,x,u). \)
\end{proof}

The following
proposition provides a situation in a semigroup context where a
much stronger version of the condition (i) in Theorem~\ref{analth} applies.
It also sheds light on the connection between semi-group theory
and holomorphic properties of  generalized transforms.

\begin{prop}
\label{SG_AV}
Let
$C_0(\mathbb{R}^n)$ be the Banach space of
continuous functions $f:$ $\mathbb{R}^{n}\rightarrow\mathbb{C}$ which vanish at infinity, equipped with supremum
norm:
$||f||:=\sup_{x\in \mathbb{R}^n}|f(x)|.$ Let $A:$ $\mathfrak{D}(A)\subset C_0(\mathbb{R}^n)%
\rightarrow C_0(\mathbb{R}^n)$ be the generator of the
Feller-Dynkin  semi-group $(P_{s})_{s\geq0}$ associated with the
process $X,$ i.e. $P_{s}f(x)=$ $\E\,[f(X_{s}^{0,x})],$ $f\in C_0(\mathbb{R}^n).$ Suppose that the family $\mathfrak{F}$ is such
 that $f_{u}\in\mathfrak{D}(A^{k})$ for each $u\in I$  and
all integer $k\geq 0,$  and that for each $u\in I,$ 
\[
\sum_{k=0}^{\infty}\frac{s^{k}}{k!}\left\Vert A^{k}f_{u}\right\Vert
<\infty,\quad0\leq s<R_{u}.
\]
Then for each $u\in I,$
\begin{equation}
P_{s}f_{u}=\sum_{k=0}^{\infty}\frac{s^{k}}{k!}A^{k}f_{u},\quad0\leq
s<R_{u},\label{HG}%
\end{equation}
with convergence in $C_0(\mathbb{R}^n).$ Thus, the map $s\rightarrow P_{s}f_{u}$
for $0\leq s<R_{u}$ extends via (\ref{HG}) to the complex disc $D_{0}%
:=\left\{  s\in\mathbb{C}:\left\vert s\right\vert <R_{u}\right\}  .$ In
particular, for each $x$\ $\in\mathbb{R}^{n}$ the map $s\rightarrow P_{s}%
f_{u}(x)$ is holomorphic in $D_{0}.$
Moreover, for each $t\geq 0,$ we may extend the map $s\rightarrow P_{t+s}f_{u},$
$0\leq s<R_{u}$ to the disc $D_{0}$ via,
\begin{equation}
P_{s+t}f_{u}=P_{t}P_{s}f_{u}=\sum_{k=0}^{\infty}\frac{s^{k}}{k!}P_{t}A^{k}%
f_{u}=\sum_{k=0}^{\infty}\frac{s^{k}}{k!}A^{k}P_{t}f_{u},\quad s\in
D_{0}.\label{HG1}%
\end{equation}
\end{prop}

\begin{proof}
See Appendix.
\end{proof}
Thus, under the conditions of Proposition \ref{SG_AV},
$\mathfrak{F}=\{f_{u},$ $u\in I\}$ is a set of analytic
vectors for the generator $A$ in the sense of Definition~\ref{AV} with $\mathfrak{X}=\mathbb{R}^n$. Moreover, due to Proposition
\ref{SG_AV} the map
\[
s\rightarrow P_{s+t}f_{u}(x)=\sum_{k=0}^{\infty}\frac{s^{k}}{k!}P_{t}%
A^{k}f_{u}(x)=\sum_{k=0}^{\infty}\frac{s^{k}}{k!}A^{k}E[f_{u}(X_{t}^{0,x})],
\]
is holomorphic in $D_{0}$ for each $x\in\mathbb{R}^{n}$ and hence Theorem~\ref{analth}-(i)
is fulfilled.
\par
In this paper we do not stick to the semigroup framework because we  want to  avoid the
narrow corset conditions of Proposition~\ref{SG_AV}. We also want to
consider operators \( A \) with unbounded (for instance, affine) coefficients and sets $\mathfrak{F}$ of functions that do
not vanish at infinity (for example, (\ref{SF})). Such situations may
lead to the violation of condition
\( f_{u}\in\mathfrak{D}(A^{k}),\, k\in \mathbb{N} \) in the sense of Proposition~\ref{SG_AV}.
In particular, in the next Sections~\ref{genA}-\ref{LAFF} we will focus
on general operators $A$ with affine coefficients   and in Section~\ref{ItoLevy}
on affine processes related to affine generators satisfying a kind of admissibility conditions.

\section{Affine generators}\label{genA}

Let us now consider generators of the form (\ref{Gen}) with affine coefficients. It is important to note that in this section $A$ may or may not be a
generator of some Feller-Dynkin process.
The next theorem and its corollaries
say that the Fourier
basis \( \mathfrak{F}=\{ f_{u}(x):=e^{iu^{\top}x},\ u\in \mathbb{R}^n \} \) is a set of analytical
vectors for an operator $A$ of the form (\ref{FBK2}), where the coefficients $\mathfrak{a}_\alpha(x)$ are affine functions of $x$ and
satisfy certain growth conditions for $|\alpha|\rightarrow\infty.$ Moreover, an explicit estimate for the radius of convergence  is given. 

\begin{thm}\label{affthm}
Let  \( A \) be a generator of the form (\ref{FBK2}) with affine coefficients $\mathfrak{a}_\alpha(x),$ i.e.
for all multi-indexes $\alpha,$
\begin{equation}
\mathfrak{a}_{\alpha}(x)=:c_{\alpha}+x^{\top}d_{\alpha}\label{a-c-d}, \quad x\in\mathfrak{X},
\end{equation}
where $\mathfrak{X}\subset\mathbb{R}^n$ is an open region, $c_{\alpha}$ is a
scalar constant, and $d_{\alpha}\in\mathbb{R}^{n}$ is a constant vector.  
Assume that the series
\begin{equation}
\sum_{\left\vert \alpha\right\vert >0}\mathfrak{a}_{\alpha}(x)\left(
\mathfrak{i}u\right)^{\alpha}
\label{momres}
\end{equation}
converges absolutely  for all $u\in\mathbb{R}^{n}.$
Then, for every $u\in \mathbb{R}^{n}$ and $x\in\mathfrak{X}$ it holds
\begin{equation}\label{coefest}
\left\vert A^{r} f_u(x)\right\vert\leq (r+1)!\, 2^{nr}(1+\| x \|)^{r}\theta^{r}(\|u\|)
\end{equation}
with $\Vert x\Vert=\displaystyle\max_{i=1,\ldots,n}|x_{i}|$,
\begin{eqnarray}
\label{theta}
\theta(v):=\sum_{k\geq1}2^k(1+v)^{k}\mathcal{D}_{k}^{\mathfrak{a}}, \quad v\in \mathbb{R}_{+}
\end{eqnarray}
and
\[
\mathcal{D}_{k}^{\mathfrak{a}}:=\sup_{x\in \mathfrak{X}}\,
\max_{|\alpha|=k,\,|\beta|\leq1}\frac{|
\partial_{x^{\beta}}\mathfrak{a}_{\alpha}(x)|}{1+\|x\|}.
\]
\end{thm}
The proof of Theorem~\ref{affthm} is given in the Appendix.

\begin{cor}
If in Theorem~\ref{affthm} the region $\mathfrak{X}$ is bounded, the Fourier basis
\( \mathfrak{F} \) constitutes a set of analytic vectors for the affine operator $A$ in  $\mathfrak{X}.$
\end{cor}
\begin{cor}
If in Theorem~\ref{affthm} there exists
for any $\varsigma>0$ a constant $M$ (which may depend on
$\varsigma>0$) such that
\[
\mathcal{D}_{k}^{\mathfrak{a}}\leq M\varsigma^{k}/k!\quad
k\geq 1,
\]
then%
\[
\theta(v)\leq 
M\exp\left( 2\varsigma(1+v)\right).
\]
\end{cor}
\begin{cor}
If in Theorem~\ref{affthm} it holds that \( \mathfrak{a}_{\alpha}(x)\equiv 0 \)
for \( |\alpha|>2 \) (generator of diffusion type)
then%
\[
\theta(v)\leq 
C(1+v)^{2}, \quad C>0.
\]
\end{cor}
\begin{rem}
The requirement that \eqref{momres} converges for all \( u \) imposes restrictions
only on the tails of the measure \( \nu \) and so does not exclude infinite activity processes.
\end{rem}
For an affine operator $A$ the sequence (\ref{qprod}) can be explicitly constructed via the next proposition,
which is proved in the Appendix.

\begin{prop}\label{affprop}
Let $A$ be an affine generator as in Theorem~\ref{affthm} and define
\begin{align*}
\mathfrak{b}_{\beta}(x,u)  &:=\partial_{u^\beta} \frac{Af_u(x)}{f_u(x)}=
\sum_{\alpha\ge0}\mathfrak{a}_{\alpha+\beta}(x)
\frac{(\alpha+\beta)!}{\alpha!}(\mathfrak{i}u)^{\alpha}\\
&=:
\mathfrak{b}_{\beta}^{0}(u)+\sum_{\kappa,\,\left\vert \kappa\right\vert
=1}\mathfrak{b}_{\beta,\kappa}^{1}(u)\,x^{\kappa},
\end{align*}
with $\mathfrak{a}_0:=0.$ We set $A^r f_u(x)$ $=:$ $g_r(x,u) f_u(x) $ and, for fixed $\eta$ $>$ $0,$ $q^{(\eta)}_r(x,u)=:h_r(x,u) f_u(x),$ where
both $g_r$ and $h_r$ are polynomials in $x$ of degree $r.$ It holds  
\begin{equation}\label{ghpol}
g_{r}(x,u)=:\sum_{\left\vert \gamma\right\vert \leq r}g_{r,\gamma}(u)x^{\gamma},\quad
h_{r}(x,u)=:\sum_{\left\vert \gamma\right\vert \leq r}h_{r,\gamma}(u)x^{\gamma},%
\end{equation}
where  $g_r$ and $h_r$ satisfy $g_0$  $\equiv$ $g_{0,0}$ $\equiv$ $h_0$ $\equiv$ $h_{0,0}$ $\equiv$ $1,$ and for $r\ge0,$ respectively,
\begin{eqnarray}
g_{r+1,\gamma}&=& \label{grecu}
\sum_{\left\vert \beta\right\vert \leq r-\left\vert \gamma\right\vert
}\binom{\gamma+\beta}{\beta}g_{r,\gamma+\beta}\mathfrak{b}_{\beta}^{0}%
\\
&+&\sum_{\left\vert \kappa\right\vert =1,\, \kappa\leq\gamma}\ \sum_{\left\vert
\beta\right\vert \leq r+1-\left\vert \gamma\right\vert }\binom{\gamma-\kappa+\beta}{\beta}
g_{r,\gamma
-\kappa+\beta}\mathfrak{b}_{\beta,\kappa
}^{1},\ {\rm and,} \notag\\
(r+1)h_{r+1,\gamma}&=& \notag
\sum_{\left\vert \beta\right\vert \leq r-\left\vert \gamma\right\vert
}\eta^{-1}\binom{\gamma+\beta}{\beta}h_{r,\gamma+\beta}\mathfrak{b}_{\beta
}^{0}  \\
  &+&\!\! \!\! \!\! \sum_{\left\vert \kappa\right\vert =1,\, \kappa\le\gamma}\ \sum_{
  \left\vert \beta\right\vert \leq r+1-\left\vert \gamma\right\vert }%
\!\! \! \eta^{-1}\binom{\gamma-\kappa+\beta}{\beta}
h_{r,\gamma-\kappa+\beta}\mathfrak{b}_{\beta,\kappa}^{1}+rh_{r,\gamma}(u) , \notag
\end{eqnarray}
where $\left\vert \gamma\right\vert \leq r+1,$ and empty sums are defined to be zero.
\end{prop}

\begin{rem}\label{trans}{\rm
Depending on the open set $\mathfrak X$ we may consider instead of (\ref{ghpol}) for an $x_0$ $\in$ $\mathfrak X$
expansions in $x-x_0$ rather than in $x.$ For simplicity we  henceforth assume $\{ 0 \}\in \mathfrak X$ which, if necessary, may be realized
by a translation of the state space. }
\end{rem}
A natural question is whether affine generators are the only ones
for which the Fourier basis constitutes a set of analytic vectors. For this paper we leave this issue as an open problem but the following
proposition shows that at any case the set of such generators is rather ``thin''.
\par
Let us put  $\mathfrak{X}=[-\pi,\pi]$ and
\begin{equation}
A=\frac{1}{2}a(x)\frac{\partial^2}{\partial x^2}+
b(x) \frac{\partial}{\partial x}.
\end{equation}
\begin{prop}\label{thin}
The set of coefficients \( (a(x),b(x)) \) such that for an arbitrary \( M>0 \)
\[
\|A^{N}f_{u} \|_{L^{2}(\mathfrak{X})}\gtrsim M^{N}N!\,,\quad N\to \infty,
\]
is dense in \( L^{2}(\mathfrak{X})\times L^{2}(\mathfrak{X}) \).
\end{prop}
\begin{proof}
Without loss of generality let us assume that \( b(x)\equiv 0 \) and $u>0.$
The general case can be considered along the same ideas and is only formally more complicated.
Any $a\in L^ 2(\mathfrak{X})$  may be approximated  (in $L^2$-sense) by a finite Fourier series
\begin{equation}
a(x)\approx\sum_{l=1}^ n a_le^ { il x}.
\end{equation}
Thus, for given $\varepsilon >0$ we can find natural \( n \) and amplitudes $a_l$ (\( a_{n}\neq 0\)) such that
\begin{equation}
\left\|a(x)-\sum_{l=1}^n a_le^ {il x}\right\|_{L^ 2(\mathfrak{X})}\leq \varepsilon .
\end{equation}
The corresponding approximative operator  is given by
\begin{equation}
\widetilde A:=\sum_{l=1}^ n \widetilde A_l= \sum_{l=1}^ n a_le^ {i l x}\frac{\partial^ 2}{\partial x^
2}.
\end{equation}
Using the fact that for any \( s_{1},\ldots,s_{k}\in \mathbb{N}, \)
\begin{equation*}
\widetilde A_{s_1}\cdots \widetilde A_{s_k}e^ {iux}=(-1)^ ka_{s_1}\cdots a_{s_k} \prod_{l=0}^ {k-1}
\left(u+\sum_{j=0}^ {l} s_{j}\right)^ 2 e^ {i(u+\sum_{j=1}^ ks_{j})x},
\end{equation*}
and setting
\begin{equation*}
F_k :=\frac{1}{2\pi}\int_{-\pi}^{\pi} e^ {-i k x}f_{u}^{-1}(x)\widetilde{A}^{N}f_{u}(x) dx, \quad k\in
\mathbb{N},
\end{equation*}
we have \( F_{k}=0 \) for \( k>nN, \) and for $N\to \infty$
\begin{equation*}
F_{nN}=(-1)^ N a_n^ {N}\prod_{l=0}^{N-1}(u+nl)^ 2
\sim (-1)^ N a_n^ {N}n^{2N}((N-1)!)^2(N-1)^{2u/n}.
\end{equation*}
Further, by Parseval's identity it holds
\[
\| \widetilde{A}^{N}f_{u} \|_{L^{2}(\mathfrak{X})}=
\left[ 2\pi\sum_{k=0}^{nN} |F_{k}|^{2} \right]^{1/2},
\]
and then we  are done.
\end{proof}
\noindent Obviously, Proposition~\ref{thin} may be formulated with respect to any compact interval.

\section{Log-affine representations for the affine Cauchy problem}\label{LAFF}
Let us consider the Cauchy problem (\ref{FBK1})
for affine generators $A$ of the form (\ref{FBK2}), under the assumption (\ref{momres}).
As in (\ref{a-c-d}) we set $\mathfrak{a}(x)$ $=$ $c_\alpha$ $+$ $x^\top d_\alpha.$ 
The ansatz
\begin{equation}
\widehat{p}(s,x,u)=\exp\left(  C(s,u)+x^{\top}D(s,u)\right)  , \label{AnsatzA}
\end{equation}
for scalar $C(s,u)$ and vector valued $D(s,u)$,
where $C(0,u)=0$ and $D(0,u)=\mathfrak{i}u,$ for the Cauchy problem (\ref{FBK1}) yields,
\begin{align*}
\partial_{s}C+x^{\top}\partial_{s}D &  =\sum_{\left\vert \alpha
\right\vert >0}\mathfrak{a}_{\alpha}(x)D^{\alpha}
  =\sum_{\left\vert \alpha\right\vert >0}c_{\alpha}D^{\alpha}+\sum
_{\left\vert \alpha\right\vert >0}x^{\top}d_{\alpha}D^{\alpha},
\end{align*}%
and so
\begin{align*}
\partial_{s}C   =\sum_{\left\vert \alpha\right\vert >0}c_{\alpha}D^{\alpha
},\quad
\partial_{s}D   =\sum_{\left\vert \alpha\right\vert >0}d_{\alpha}D^{\alpha
},\\
C(0,u)   =0,\quad D(0,u)=\mathfrak{i}u.
\end{align*}
We thus have a system of ordinary differential equations (ODEs), which reads component-wise%
\begin{align}
\partial_{s}C   =\sum_{\left\vert \alpha\right\vert >0}c_{\alpha}D^{\alpha
},\quad
\partial_{s} D_{j}  =\sum_{\left\vert \alpha\right\vert >0}d_{\alpha}%
^{(j)}D^{\alpha},\quad j=1,...,n,\nonumber\\
C(0,u)   =0,\quad D_{j}(0,u)=\mathfrak{i}u_{j}\label{Ricatti}.
\end{align}
By assumption (\ref{momres}), the series
\begin{equation}\label{cda}
\sum_{\left\vert \alpha\right\vert >0}c_{\alpha}y^{\alpha},\quad
\sum_{\left\vert \alpha\right\vert >0}d_{\alpha}^{(j)}y^{\alpha},\quad
j=1,...,n,
\end{equation}
are absolutely convergent for all $y\in\mathbb{R}^{n},$ and thus define
terms-wise differentiable $C^{(\infty)}(\mathbb{R}^{n})$ functions. In
particular, they are locally Lipschitz and so according to standard ODE theory
the system (\ref{Ricatti}) has for fixed $u\in\mathbb{R}^{n}$ a unique solution
$(C(s,u),D(s,u))$ for $0\leq s<s_{u}^{\infty}\leq\infty,$ where
$(s,C(s,u),D(s,u))$ leaves any compact subset of $\mathbb{R}\times
\mathbb{R}\times\mathbb{R}^{n},$ when $s\uparrow s_{u}^{\infty}.$

\begin{rem}
By a general theorem from analysis (e.g., see \cite{D}), it follows that the solution of
(\ref{Ricatti}) extends component-wise  holomorphically in $s$ into a
disc around $s=0,$ due to the analyticity of (\ref{cda}). This implies that (\ref{AnsatzA})
is holomorphic  in  $s.$  So, besides Theorem~\ref{affthm}, also
along this line one may show  that (\ref{AnsatzA}) can be represented as a power series
 of the form (\ref{ConvP}). I.e.,   in particular,
the Fourier basis (\ref{SF}) constitutes a set of analytic vectors for the affine operator $A.$
However, the direct approach in the proof  of Theorem~\ref{affthm} (see Appendix) leads
to an explicit estimate \eqref{coefest} and
allows for investigating possible extensions of \( \widehat p(s,x,u) \) into
a strip   containing the real axis in  the complex plane (see Theorem~\ref{affineD}).
Moreover, it also suggests the line to follow in cases where $A$
is not affine and/or the function base is not of the form (\ref{SF}).
\end{rem}

 Let us suppose that for fixed $u\in\mathbb{R}^{n}$  the statements of  Theorem~\ref{analth} hold. Then  we  obtain for $0\le s<s^\infty_u\le\infty,$
$x\in\mathfrak{X,}$
\begin{align}
\widehat{p}(s,x,u)&=\exp\left(  C(s,u)+x^{\top}D(s,u)\right) \label{logC}\\
 &=\sum
_{k=0}^{\infty}q_{k}^{(\eta_{u})}(x,u)(1-e^{-\eta_{u}s})^{k}.
\notag
\end{align}
 {\rm  Since $q_{0}^{(\eta_{u})}(x,u)=f_{u}(x)=\exp\left(
\mathfrak{i}u^{\top}x\right)  \neq0$ we have, taking into account the
boundary conditions for $C$ and $D,$ at least for small enough $\varepsilon>0,$
\begin{align}\label{epsi}
C(s,u)+x^{\top}D(s,u)  & =\sum_{k=0}^{\infty}\rho_{k}%
^{(\eta_{u})}(x,u)(1-e^{-\eta_{u}s})^{k}, \quad 0\le s<\varepsilon,
\end{align}
where  by a standard lemma on the power series expansion of the logarithm of
a power series,
\begin{align}
\rho_{0}^{(\eta_{u})}(x,u)  & =\ln q_{0}^{(\eta_{u})}(x,u)=\mathfrak{i}%
u^{\top}x,\label{rhoaff}\\
\rho_{k}^{(\eta_{u})}(x,u)  & =\frac{1}{f_{u}(x)}\left[  q_{k}^{(\eta_{u}%
)}(x,u)-\frac{1}{k}\sum_{j=1}^{k-1}j\rho_{j}^{(\eta_{u})}(x,u)q_{k-j}%
^{(\eta_{u})}(x,u)\right],  \quad k\geq1.\notag
\end{align}
Thus by (\ref{epsi}),   the $\rho_{k}^{(\eta_{u})}$
are necessarily affine in $x$!
}

\begin{rem}
{\rm  It is possible to prove directly  that the functions $\rho_{k}^{(\eta
_{u})}$  defined above
 are affine in $x$ using Proposition~\ref{affprop} via  a (rather laborious) induction procedure, so
without using a local solution of (\ref{Ricatti}).
}
\end{rem}

\begin{thm}\label{logbran}
Suppose that for fixed $u\in I$ the statement of Theorem~\ref{analth}-(i)
holds  for an open region  $\mathfrak{X}$ and, in addition,
for $s\in G_{R_u}$ and $x\in\mathfrak{X}$ it holds $\widehat p(s,x,u)\neq0.$
Then, with
$\rho_{k}^{(\eta_{u})}(x,u)=:\rho_{k}^{(\eta_{u},0)}(u)+x^\top\rho_{k}^{(\eta_{u},1)}(u)$ determined  by (\ref{rhoaff}),
 we have $$\widehat p(s,x,u)=\exp\left[\sum_{k=0}^{\infty}\left(\rho_{k}^{(\eta_{u},0)}
 +x^\top \rho_{k}^{(\eta_{u},1)}(u)\right)(1-e^{-\eta_{u}s})^k\right],\quad 0\le s<\infty.$$
\end{thm}

\begin{proof} Let $u\in I$ and $x\in\mathfrak{X}$ be fixed.
Since $G_{R_u}$ is simply connected and $s\rightarrow\widehat p(s,x,u)$ is holomorphic  and non-zero in  $G_{R_u}$, there exists a branch $s\rightarrow L(s,x,u)$ of the logarithm such that $\widehat p(s,x,u)=\exp(L(s,x,u))$ for all $s\in G_{R_u}.$ Along the same line  as in Theorem~\ref{analth} we then argue that  there exists an $\eta_u>0$ such that
$w\rightarrow L(\Phi_{\eta_u}(w),x,u) $ (see the proof of Theorem~\ref{analth}) is holomorphic in the unit disc $\{w:\mid w\mid<1\},$ hence, there exists $\widetilde\rho_{k}(x,u)$ such that
$$
L(\Phi_{\eta_u}(w),x,u)=\sum_{k\ge0}\widetilde\rho_{k}(x,u)w^k, \quad 0\le |w|<1,\quad {\rm and\ so}
$$
$$
L(s,x,u)=\sum_{k\ge0}\widetilde\rho_{k}(x,u)(1-e^{-\eta_{u}s})^k,\quad 0\le s<\infty.
$$
Since the later expansion must coincide with (\ref{epsi}) for small $s,$ it follows that necessarily $\widetilde\rho_{k}(x,u)$ $=$ $\rho^{(\eta_u)}_{k}(x,u)$
and the theorem is proved.
 \end{proof}

Let us now pass to another interesting log-affine representation for the characteristic function.
From (\ref{ghpol}) and (\ref{grecu}) we derive formally
\begin{align*}
 \sum_{r\geq0}q_{r}^{(\eta)}(x,u)(1-e^{-\eta s}%
)^{r}&=e^{\mathfrak{i}u^{\top}x}\sum_{r\geq0}\sum_{\gamma\geq
0}h_{r,\gamma}(u)x^{\gamma}(1-e^{-\eta s})^{r}1_{\left\vert \gamma\right\vert
\leq r}\\
& =e^{\mathfrak{i}u^{\top}x}\sum_{\gamma\geq0}x^{\gamma}\sum_{r\geq
0}h_{\left\vert \gamma\right\vert +r,\gamma}(u)(1-e^{-\eta s})^{\left\vert
\gamma\right\vert +r}.%
\end{align*}
Suppose  that the requirements of Theorem~\ref{logbran} hold.  Then, using Theorem~\ref{affthm}, it is easy to show that for small enough $\varepsilon>0$,
$$
\sum_{\gamma\geq0}\|x\|_\infty^{|\gamma|}\sum_{r\geq
0}|h_{\left\vert \gamma\right\vert +r,\gamma}(u)||w|^{\left\vert
\gamma\right\vert +r}<\infty, \quad {\rm if}\ \quad |w|<\varepsilon,\quad \|x\|_\infty<\varepsilon
$$
(see Remark~\ref{trans}).
Thus, for $|s|$ and $\|x\|_\infty$ small enough  we obtain
\begin{eqnarray*}
\ln\widehat{p}(s,x,u)&=&\mathfrak{i}u^{\top%
}x+\ln\left(\sum_{\gamma\geq0}x^{\gamma}\sum_{r\geq0}h_{\left\vert \gamma\right\vert
+r,\gamma}(u)(1-e^{-\eta s})^{\left\vert \gamma\right\vert +r}\right)\\
&=&C(s,u)+x^{\top}D(s,u),%
\end{eqnarray*}
with (in multi-index notation)
\begin{align}
C(s,u)  & =\ln\left(\sum_{r\geq0}h_{r,0}(u)(1-e^{-\eta_{u}s})^{r}\right),\label{Comp}\\
D^{\kappa}(s,u)  & =\mathfrak{i}u^{\kappa}+\frac{\sum
_{r\geq1}h_{r,\kappa}(u)(1-e^{-\eta_{u}s})^{r}}{\sum_{r\geq0}h_{r,0}%
(u)(1-e^{-\eta_{u}s})^{r}},\quad\left\vert \kappa\right\vert =1.\notag
\end{align}
However, the left- and right-hand-sides of (\ref{Comp}) are holomorphic for all $s$ $\in$ $G_{R_u}$ and we  so arrive
at the representation
\begin{align}
\widehat{p}(s,x,u)  & =\exp\left[  \ln\left(\sum_{r\geq0}h_{r,0}(u)(1-e^{-\eta_{u}%
s})^{r}\right)+\mathfrak{i}u^{\top}x\right. \label{logseries} \\
& +\left. x^{\top}  \frac{\sum_{r\geq1}h_{r}%
(u)\,(1-e^{-\eta_{u}s})^{r}}{\sum_{r\geq0}h_{r,0}(u)(1-e^{-\eta_{u}s})^{r}%
}\right],\quad s\in G_{R_u,}\ x\in\mathfrak{X,}\notag
\end{align}
with%
\[
h_{r}(u):=\left[  h_{r,i}(u)\right]  _{i=1,...,n},%
\]
where for $1\leq i\leq n,$ the multi-index $(\delta_{ij})_{j=1,...,n}$ is
identified with $i.$
\ \\

Particularly due to the explicit estimate (\ref{coefest}) for affine generators in Theorem~\ref{affthm}, we may proof the next theorem
(which is a non-probabilistic version of Proposition~\ref{probcrit}   in the situation where $A$ is affine).
\begin{thm}
\label{affineD}
Let \( \mathfrak{X} \) be a bounded domain. Assume that the system {\rm (\ref{Ricatti})} is
non-exploding, i.e. $s_{u}^{\infty}$ $=$ $\infty,$ and that for any
fixed $u\in\mathbb{R}^n$ the solution
$D(s,u)$ remains bounded as $s$ $\rightarrow$ $\infty.$ Then, there exists a radius
$R_u>0$ such that for any $t\ge 0$ the map $s$ $\rightarrow$ $\widehat{p}(t+s,x,u),$
$0\le s<R_u$ has a holomorphic extension to the disc $\{s\in \mathbb{C}:\vert s\vert<R_u\}.$ Moreover, it holds
\begin{equation}
 \widehat p({t+s},x,u)=\sum_{k=0}^{\infty}\frac{s^{k}}%
{k!}A^{k}\widehat p({t},\cdot,u)\,(x),\quad |s|<R_{u},\quad x\in \mathfrak{X} .
\label{ConvP1'}%
\end{equation}
\end{thm}
\begin{rem}
\label{RLB}
The maximal extension radius \( R_{u} \) satisfies
\begin{equation}
\label{RU}
R_u\geq  \frac{1}{2^{n}\,\theta(\|D^{*}(u)\|)}\inf_{x\in \mathfrak{X}}\frac{1}{1+\|x\|},
\end{equation}
where function \( \theta\) is defined in \eqref{theta} and
\( \|D^{*}(u)\|=\sup_{s>0}\|D(s,u)\| \).
\end{rem}
\begin{proof}
Denote the right-hand-side of \eqref{ConvP1'} by \( {\widetilde p}(t,s,x,u) \).
Obviously, ${\widetilde p}(t,0,x,u)={\widehat p}(t,0,x,u).$ Let us define
\[
{\widetilde p}^{(N)}(t,s,x,u):=\sum_{k=0}^{N}\frac{s^{k}}{k!}A^{k}\widehat
p({t},\cdot,u)\,(x).
\]
Since  $\widehat p({t},x,u)$ $=$ $\exp(C(t,u)+x^\top D(t,u)),$  Theorem~\ref{affthm} implies
that the series
\begin{equation}
\sum_{k=0}^{\infty}\frac{s^{k}}{k!}A^{k}\widehat p({t},\cdot,u)\,(x)
\end{equation}
is absolutely and uniformly convergent on any compact subset of $\mathfrak{X}\times\{s\in \mathbb{C}:\vert s\vert<R_u\}$
if \( R_{u} \) satisfies \eqref{RU}.
So, both ${\widetilde p}^{(N)}(t,s,x,u)$ and
\[
A{\widetilde p}^{(N)}(t,s,x,u)=\sum_{k=0}^{N}\frac{s^{k}}{k!}A^{k+1}\widehat p({t},\cdot,u)\,(x)
\]
converge for $N\rightarrow\infty$ uniformly over any compact subset of
$\mathfrak{X}$ and $s$ with $\left\vert s\right\vert <R_{u}-\varepsilon$ for
any $\varepsilon>0$.
Hence, for $\left\vert s\right\vert <R_{u}-\varepsilon$
\begin{align*}
\frac{\partial}{\partial
s}{\widetilde p}^{(N)}(t,s,x,u)&=\sum_{k=0}^{N-1}\frac{s^{k}}%
{k!}A^{k+1}\widehat p(t,\cdot,u)=A{\widetilde p}^{(N-1)}(t,s,x,u),\\
{\widetilde p}(t,0,x,u)  & =\widehat{p}(t,x,u),
\end{align*}
and thus, by closeness of the operator $A$ and  uniqueness of the
Cauchy problem (\ref{FBK1})-(\ref{FBK2}), we have \( {\widetilde p}(t,s,x,u)=\widehat{p}(t+s,x,u). \)
\end{proof}

\section{Full expansion of a specially structured one dimensional affine
system}\label{EX}

\bigskip Let us consider Cauchy problem \eqref{FBK1} for $n=1$ with $f_{u}%
(x)=\exp(\mathfrak{i}ux),$ where the jump-kernel in the generator $A$ (see
\eqref{Gen}) has a special affine structure of the form%
\[
v(x,dz)=:(\lambda_{0}+\lambda_{1}x)\mu(dz),
\]
and where the diffusion coefficients have a similar structure,%
\[
b(x)=(\lambda_{0}+\lambda_{1}x)\theta,\quad a(x)=(\lambda_{0}+\lambda
_{1}x)\vartheta,
\]
for some constants $\lambda_{0},\lambda_{1},\theta,\vartheta\in\mathbb{R},$ and
measure $\mu.$ So, in Proposition~\ref{affprop} the $\mathfrak{a}_{\alpha}$
have the form $\mathfrak{a}_{l}$ $=:$ $(\lambda_{0}+\lambda_{1}x)\eta_{l}$
where%
\[
\eta_{0}:=0,\quad\eta_{1}=\theta,\quad\eta_{2}:=\frac{1}{2}\left(
\vartheta+\int z^{2}\mu(dz)\right)  ,\quad\eta_{l}:=\frac{1}{l!}\int z^{l}%
\mu(dz),\quad l>2.
\]
Hence, in Proposition~\ref{affprop}, the $\mathfrak{b}_{\beta}$ in \eqref{grecu}
have the form%
\begin{align*}
\mathfrak{b}_{r}(x,u) &  =\mathfrak{b}_{r}^{0}(u)+x\mathfrak{b}_{r}^{1}(u)\\
&  =:(\lambda_{0}+\lambda_{1}x)\sum_{l\geq0}\eta_{l+r}\frac{(l+r)!}%
{l!}(\mathfrak{i}u)^{l}=:(\lambda_{0}+\lambda_{1}x)\frac{d^{r}}{du^{r}}\mathfrak{h}(u)
\quad r\geq 0,
\end{align*}
where $\mathfrak{h}(u):=\sum_{l\geq0}\eta_{l}(\mathfrak{i}u)^{l}.$ It is now
possible to show via \eqref{grecu} that for $r\geq1,$%
\begin{multline}
\label{GR}
g_{r}(x,u)=\sum_{\substack{p>0,\,q\geq0\\ 0<n_{1}<...<n_{q},\,\,m_{1}%
,...,m_{q}\geq0,\\r=p+n_{1}m_{1}+\cdot\cdot\cdot+n_{q}m_{q}}}\frac{1}{r!}%
\pi_{(n_{1},m_{1}),\cdot\cdot\cdot,(n_{q},m_{q})}^{(p)}\lambda_{1}%
^{r-p}(\lambda_{0}+\lambda_{1}x)^{p}\\
\times\ \mathfrak{h}^{p}(u)\
\prod\limits_{j=1}^{q}
\left(  \mathfrak{h}^{n_{j}-1}(u)\frac{d^{n_{j}}}{du^{n_{j}}}\mathfrak{h}%
(u)\right)  ^{m_{j}},
\end{multline}
with the following integer recursion procedure:
\begin{description}
\item[Initialization:]
$\pi^{(p)}\equiv1,\quad\pi_{(n_{1},m_{1}),\cdot\cdot
\cdot,(n_{q},m_{q})}^{(0)}\equiv0,\quad p,q\geq1.$
\end{description}
For all \( n_{i}>0, \, m_{i}\geq 0,\) with \( 1\leq i\leq q,  \, p,q\geq 1  \):
\begin{description}
\item[Reduction rule I:] If $m_{j}=0,$ for some \( j \), $1\leq j\leq q,$ then
\[
\pi_{(n_{1},m_{1}),\cdot\cdot\cdot,(n_{q},m_{q})}^{(p)}=\pi_{(n_{1}%
,m_{1}),\cdot\cdot\cdot,(n_{j-1},m_{j-1}),(n_{j+1},m_{j+1}),\cdot\cdot
\cdot,(n_{q},m_{q})}^{(p)}.
\]%
\item[Reduction rule II:]
\begin{multline*}
\qquad\pi_{(n_{1},m_{1}),\cdot\cdot\cdot,(n_{q-1}%
,m_{q-1}),(n_{q},m_{q})}^{(p)}=
\\
\sum_{j=1}^{q}\binom{p+n_{j}-1}{n_{j}}%
\pi_{(n_{1},m_{1}),\cdot\cdot\cdot,(n_{j},m_{j}-1),\cdot\cdot\cdot
,(n_{q},m_{q})}^{(p+n_{j}-1)}
+\pi_{(n_{1},m_{1}),\cdot\cdot\cdot,(n_{q},m_{q})}^{(p-1)}.%
\end{multline*}
\end{description}
In fact, the above recursion procedure follows automatically after substituting
\eqref{GR} as  ansatz into \eqref{grecu}.

Note that
\begin{align*}
\mathfrak{h}(u) &  =\mathfrak{i}u\theta-\frac{1}{2}\vartheta u^{2}+\sum
_{l\geq2}\frac{(\mathfrak{i}u)^{l}}{l!}\int z^{l}\mu(dz)\\
&  =\phi(u)-1+\left(  \mathfrak{i}\theta-\phi^{\prime}(0)\right)  u-\frac
{1}{2}\vartheta u^{2},
\end{align*}
where $\phi$ is the characteristic function of the measure $\mu.$ Hence,
$\mathfrak{h}$ and all its derivatives may be computed from $\phi.$
Subsequently we obtain the $q_{k}^{(\eta)}$ for the series expansion in
Theorem~\ref{analth} by \eqref{qseq}, i.e.%
\[
q_{k}^{(\eta)}(x,u)=:\frac{1}{k!}\sum_{r=0}^{k}c_{k,r}\eta^{-r}g_{r}
\]

\section{Application to affine processes}\label{ItoLevy}
Affine  processes have become very popular in recent years due
to their analytical tractability in the context of option pricing,
and their rather rich  dynamics.
Many well-known models such as Heston and Bates models fall into the class
of affine jump diffusions. Option pricing in these models is usually done
via the Fourier method which requires knowledge of the
Fourier transform of the process in closed form (see e.g. \cite{DPS}).
The functional series representations for affine generators developed
in this paper, in particular (\ref{logseries}), can be directly applied
to affine processes.
Let us recall the characterization of a regular affine process as given in
\cite{DFS}.

\begin{defi}
We  call a strong Markov process \( \{ X_{t} \} \) with generator \( A \) a regular
affine process if \( A \) is of the form \eqref{Gen}  and all functions
\[
 a_{ij}(x),\ b_{i}(x),\ v(x,dz)\quad i,j=1,\ldots,m
\]
are affine in \( x \) (see \eqref{a-c-d}),
and satisfy the set of admissability conditions spelled out in
\cite[Definition 2.6]{DFS}. These conditions
guarantee that $A$ is the generator of a
Feller-Dynkin (strong) Markov process $X$ in a subspace
of the form \( \mathbb{R}^l\times\mathbb{R}^{n-l}_{+}\subset \mathbb{R}^{n} \)
for some \( 0\leq l\leq n \).
\end{defi}

The next theorem provides a sufficient condition for convergence of
the series representation in Theorem~\ref{analth}-(ii), hence
representation \eqref{logseries}, for regular affine processes.

\begin{thm}
Let \( \{ X_s \} \) be a regular affine  process   which has a
non-degenerated limiting distribution for $s\rightarrow\infty$, and
has a generator \( A \)
which satisfies  the moment condition (\ref{momres}). Then the (conditional)
characteristic function
$\widehat p(s,x,u)$ $=$ $\E[f_u(X^{0,x}_s)],$ with $f_u(x)$
$=$ $e^{\mathfrak i u^\top x},$ has a representation according to
Theorem~\ref{analth}-(ii):
\[
\widehat{p}(s,x,u)=\sum_{k=0}^{\infty}q_{k}^{(\eta_{u})}(x,u)\left(
1-e^{-\eta_{u}s}\right)  ^{k},\quad 0\le s<\infty.%
\]
Moreover, the scaling factor \( \eta_{u} \) may be chosen according to the inequality:
\[
\eta_{u}\geq C\, \theta(L(1+\|u\|^{2}_{2})),
\]
where the (monotonic) function \( \theta\) is defined in \eqref{theta},
\(L>0\) is a constant independent of $x$, and \( C \) is a constant generally depending on \( x \).
\end{thm}
\begin{proof}
Following \cite{DFS},  $\widehat p(s,x,u)$ has  representation of the form (\ref{AnsatzA})
for $0\le s<\infty.$ The existence of a limiting distribution  implies in particular
that $D(s,u)$ in (\ref{AnsatzA}) is bounded for all $s\ge0.$ Moreover,
as shown in \cite[Section 7]{DFS}, $\widehat p(s,x,u)$ is the characteristic
function of some infinitely divisible
distribution for all \( s>0 \) (hence also in the limit \( s\to \infty \)).
As a consequence (see \cite{SA}), there exists an \( M>0 \) independent of \( s \)
such that
\[
\lim_{\| u \|_2\to \infty}\| u \|^{-2}_{2}\left|\log \widehat p(s,x,u) \right|<M.
\]
This implies that \(\| D(s,u) \|\leq L (1+\| u \|^{2}_{2})\) for some constant \( L>0 \) not depending
on \( x \) and \( s\geq 0 \). Now we apply Theorem \ref{affineD} and Remark \ref{RLB}.
\end{proof}
\begin{rem}
The existence of a limiting stationary distribution is a sufficient condition for the
boundedness of \( D(s,u) \). In fact, there are  affine processes
which have  no limit distribution but bounded \( D(s,u) \) (a trivial example
is standard Brownian motion). The study of
existence of  limiting (and stationary) distributions
for affine processes is currently under active research,
e.g. see  \cite{KS} or \cite{KU}.
\end{rem}

\section{Appendix}

\subsection*{Proof of Theorem \ref{analth}}

$(i)\Longrightarrow(ii):$  Let $\mathcal{U}$ $:=$ $\left\{  z\in\mathbb{C}:\left\vert z\right\vert
<1\right\}$ be the unit disc in the complex plane. Consider for $\eta>0$ the map%
\[
\Phi_{\eta}:z\longrightarrow-\frac{1}{\eta}\text{Ln}(1-z).
\]
Obviously, there exists an $\eta_{u}>0$ such that%
\[
(0-,\infty)\subset\Phi_{\eta_{u}}(\mathcal{U})\subset G_{R_{u}}.%
\]
Moreover, the map $\Phi_{\eta_{u}}$ is injective on $\mathcal{U}.$ Thus,
(denoting the extension in $(i)$ with $\widehat{p}$ as well) for
each $x\in\mathfrak{X},$ the function $\widehat{p}(\Phi_{\eta_{u}}(w),x,u)$ is
holomorphic in $\mathcal{U}$ and has a series expansion
\[
w\longrightarrow\widehat{p}(\Phi_{\eta_u}(w),x,u)=:\sum_{k=0}^{\infty}\widetilde{q}_{k}(x,u)w^{k},\quad\left\vert w\right\vert <1\mathfrak{,}%
\]
and as a consequence,
\begin{equation}
\widehat{p}(z,x,u)=\sum_{k=0}^{\infty}\widetilde{q}_{k}(x,u)(1-e^{-\eta
_uz})^{k},\quad z\in\Phi_{\eta_u}(\mathcal{U}),\quad x\in\mathfrak{X.}%
\label{cons}%
\end{equation}
Since (\ref{cons}) holds in particular for $z\in(0-,\infty),$ we have
in a (possibly small) $\varepsilon$-disk around $z=0,$
\begin{equation}
\widehat{p}(z,x,u)=\sum_{k=0}^{\infty}\widetilde{q}_{k}(x,u)(1-e^{-\eta
_uz})^{k}=\sum_{k=0}^{\infty}\frac{z^{k}}{k!}A^{k}f_{u}(x),\quad
0\leq\left\vert z\right\vert <\varepsilon,\label{qA}
\end{equation}
due
to Proposition~\ref{propa}.
By taking $z=0$ we have
\[
\widehat{p}(0,x,u)=\widetilde{q}_{0}(x,u)=f_{u}(x).
\]
Taking derivatives at \( z=0 \) yields
\begin{align*}
\left.  \frac{\partial^{k}}{\partial z^{k}}\sum_{l=0}^{\infty}\widetilde{q}_{l}%
(x,u)(1-e^{-\eta_u z})^{l}\right\vert _{z=0}  &  =\left.  \sum_{l=0}^{k}%
\widetilde{q}_{l}(x,u)\frac{\partial^{k}}{\partial z^{k}}\sum_{j=0}^{l}\binom{l}%
{j}(-1)^{j}e^{-j\eta_u z}\right\vert _{z=0}\\
&  =\left.  \sum_{l=0}^{k}\widetilde{q}_{l}(x,u)\sum_{j=0}^{l}\binom{l}{j}(-1)^{j}%
e^{-j\eta_u z}(-j\eta_u)^{k}\right\vert _{z=0},%
\end{align*}
hence
\begin{equation}
\sum_{l=0}^{k}\widetilde{q}_{l}(x,u)\sum_{j=0}^{l}\binom{l}{j}(-1)^{j}(-j)^{k}=\eta^{-k}
_uA^{k}f_{u}(x).\label{qsys}
\end{equation}
In Lemma~\ref{TeLe} it is proved that the solution of system (\ref{qsys}) is given
by $\widetilde{q}_{l}(x,u)={q}_{l}^{(\eta_u)}(x,u)$ with ${q}_{l}^{(\eta_u)}$ defined in  (\ref{qseq}).\\
Next we prove the uniform convergence as stated in (ii). A well known property of Stirling numbers
implies that the series
\begin{align*}
\sum_{r=1}^{\infty}\left\vert \eta_u^{-r}A^{r}f_{u}(x)\right\vert \sum
_{k=r}^{\infty}c_{k,r}\frac{\left\vert w\right\vert ^{k}}{k!}  &  =\sum_{r=1}^{\infty}\left\vert
\frac{A^{r}f_{u}(x)}{r!}\right\vert \eta_u^{-r} \left\vert \log(1-|w|)\right\vert ^{r}
\end{align*}
converges uniformly on any compact subset of $\mathfrak{X}\times\{w:\left\vert w\right\vert < 1-e^{-\eta_u R_{u}}\}.$
Thus,  by a Fubini argument the series
\[
\sum_{k=0}^{\infty}\frac{(1-e^{-\eta_u s})^{k}}{k!}\sum_{r=1}^{k}c_{k,r}\eta_u^{-r}A^{r}%
f_{u}(x)=\sum_{k=0}^{\infty}q_k^{({\eta}_{u})}(x,u)(1-e^{-\eta_u s})^{k}%
\]
is also uniformly convergent on any compact subset of $\mathfrak{X} $ $\times$ $\left\{  s: 0\le s< R_u\right\}.$
\\
$(ii)\Longrightarrow(iii):$ Is obvious, take $\varepsilon_u:=1-e^{-\eta_u R_{u}}.$ \\
\\
$(iii)\Longrightarrow(i)$ Let $\eta_u$ and $\varepsilon_u$ be such that $(iii)$ holds. We may then define (see the proof of (ii))
\begin{equation}
\widetilde{p}(z,x,u)=\sum_{k=0}^{\infty}{q}^{(\eta_u)}_{k}(x,u)(1-e^{-\eta
_uz})^{k},\quad x\in\mathfrak{X},%
\label{cons1}%
\end{equation}
which is holomorphic in $z\in\Phi_{\eta_u}(\mathcal{U}).$
We first note that $\widetilde{p}(0,x,u)=f_{u}(x).$ Next we consider
for $0\leq s<-\eta_u^{-1}\ln\varepsilon_u,$%
\[
\widetilde{p}^{(N)}(s,x,u)=\sum_{k=0}^{N}q_{k}^{(\eta_{u})}(x,u)(1-e^{-\eta
_{u}s})^{k},
\]
which satisfies
\begin{align}
\frac{\partial}{\partial s}\widetilde{p}^{(N)}(s,x,u)  & =\sum_{k=1}^{N}%
kq_{k}^{(\eta_{u})}(x,u)(1-e^{-\eta_{u}s})^{k-1}\eta_{u}e^{-\eta_{u}%
s}\nonumber\\
& =-\eta_{u}\sum_{k=1}^{N}kq_{k}^{(\eta_{u})}(x,u)(1-e^{-\eta_{u}s})^{k}\nonumber\\
&+\eta_{u}\sum_{k=0}^{N-1}(k+1)q_{k+1}^{(\eta_{u})}(x,u)(1-e^{-\eta_{u}s}%
)^{k}\label{van0}\\
& =-\eta_{u}Nq_{N}^{(\eta_{u})}(x,u)(1-e^{-\eta_{u}s})^{N}\notag\\
&+\sum_{k=0}%
^{N-1}Aq_{k}^{(\eta_{u})}(x,u)(1-e^{-\eta_{u}s})^{k}\label{van}%
\end{align}
by some rearranging and using (\ref{rec}). Since due to $(iii)$ the first term in
(\ref{van}) vanishes for $N\rightarrow\infty,$ we obtain
\[
\frac{\partial}{\partial s}\widetilde{p}(s,x,u)=\lim_{N\rightarrow\infty}\frac{\partial
}{\partial s}\widetilde{p}^{(N)}(s,x,u)=\lim_{N\rightarrow\infty}%
A\widetilde{p}^{(N-1)}(s,x,u),
\]
together with%
\[
\lim_{N\rightarrow\infty}\widetilde{p}^{(N)}(s,x,u)=\widetilde{p}(s,x,u).
\]
From the uniform convergence as stated in (iii) it follows easily that the two series in (\ref{van0}), the first term in (\ref{van}), and so also the second term in (\ref{van}) converge uniformly in the same sense.  Thus, the above   limits are uniform on compacta accordingly. %
Since the operator $A$ is closed, we so obtain %
\[
\frac{\partial}{\partial s}\widetilde{p}(s,x,u)=A\widetilde{p}(s,x,u),\quad
0\leq s<-\eta_u^{-1}\ln\varepsilon_u,
\]
and  by uniqueness of the Cauchy problem associated with the operator~$A$ we thus have
\[
\widehat{p}(s,x,u)=\widetilde{p}(s,x,u)=\sum_{k=0}^{N}q_{k}^{(\eta_{u}%
)}(x,u)(1-e^{-\eta_{u}s})^{k},\quad0\leq s<-\eta_u^{-1}\ln\varepsilon_u.
\]
Because of the assumption that $\widehat{p}(s,x,u)$ is holomorphically extendable in each $s,$ $0\le s<\infty,$ we then must have $\widehat{p}(s,x,u)$ $=$ $\widetilde{p}(s,x,u)$ for $0\le s<\infty.$  Finally, it is not difficult to see that there exists $R'_u>0$ such that $G_{R'_u}$ $\subset$ $\Phi_{\eta_u}(\mathcal{U})$, hence
$(i)$ is proved. $\Box$

\begin{lem}\label{TeLe}
The solution of
\begin{equation}\label{indp}
\sum_{l=0}^{k}{q}_{l}(x,u)\sum_{j=0}^{l}\binom{l}{j}(-1)^{j}(-j)^{k}=B^{k}f_{u}(x)
\end{equation}
satisfies (\ref{qprod}), and (\ref{rec}) respectively, where $B:=\eta^{-1}A.$
\end{lem}

\begin{proof}
Suppose that the solution $q_l$ of (\ref{indp})
satisfies $(l+1)q_{l+1}$ $=$ $Bq_l+lq_l,$ see (\ref{rec}),
for $0\leq l\leq k,$ $k>0.$ Then
(note that summations may be started at $l,j=1$ for $k>0$),%
\begin{equation}
\sum_{l=1}^{k+1}q_{l}(x,u)\sum_{j=1}^{l}\binom{l}{j}(-1)^{j}(-j)^{k+1}=B^{k+1}f_{u}(x) \label{COEFK1}
\end{equation}
transforms to
\begin{gather*}
q_{k+1}(x,u)\sum_{j=1}^{k+1}\binom{k+1}{j}(-1)^{j}(-j)^{k+1}+\\
\sum_{l=1}^{k}\sum_{j=1}^{l}\binom{l}{j}(-1)^{j}(-j)^{k+1}\frac{1}{l!}%
{\displaystyle\prod\limits_{r=0}^{l-1}}
\left(B+rI\right)f_{u}(x)=B^{k+1}f_{u}(x),
\end{gather*}
and after some straightforward algebra to
\begin{gather*}
q_{k+1}(x,u)(-1)^{k+1}\sum_{j=1}^{k+1}\binom{k+1}{j}(-1)^{j}j^{k+1}=\\
B^{k+1}f_{u}(x)-\sum_{l=1}^{k}\sum_{j=1}%
^{l}\binom{l}{j}(-1)^{j}(-j)^{k+1}\frac{1}{l!}%
\left(B+rI\right)f_{u}(x).
\end{gather*}
{\em Claim:}
For any \( k>0 \)
\begin{equation}\label{COEFL1}
(-1)^{k}\sum_{j=0}^{k}\binom{k}{j}(-1)^{j}j^{k}=k!.
\end{equation}
This claim is easily proved by considering $h(s):=(1-e^{-s})^{k}$. On the one hand,
\[
\frac{d}{ds^{k}}h(s)=\frac{d}{ds^{k}}\sum_{j=0}^{k}(-1)^{j}e^{-js}=\sum
_{j=0}^{k}\binom{k}{j}(-1)^{j}(-j)^{k}e^{-js},
\]
and on the other, for $\left\vert s\right\vert $ small enough, $h(s)=(s-\frac
{1}{2}s^{2}+...)^{k}=s^{k}+....$ Therefore,%
\[
\left.  \frac{d}{ds^{k}}h(s)\right\vert _{s=0}=k!=\sum_{j=0}^{k}\binom{k}%
{j}(-1)^{j}(-j)^{k}.
\]
Using (\ref{COEFL1}),  \eqref{COEFK1} is equivalent to
\begin{gather*}
q_{k+1}(x,u)(k+1)!=B^{k+1}f_{u}(x)-\\
(-1)^{k+1}\sum_{l=1}^{k}\sum_{j=1}^{l}\binom{l}{j}(-1)^{j}j^{k+1}\frac{1}{l!}%
{\displaystyle\prod\limits_{r=0}^{l-1}}
\left(B+rI\right)f_{u}(x)
\end{gather*}
The rest follows from the next claim.\\
{\em Claim:}
For any $k>0$ and \( x\in \mathbb{R} \) it holds%
\begin{equation}
x^{k+1}-(-1)^{k+1}\sum_{l=1}^{k}\sum_{j=1}^{l}\binom{l}{j}(-1)^{j}%
j^{k+1}\left[  \frac{1}{l!}\prod_{r=0}^{l-1}(x+r)\right]={\displaystyle\prod_{r=0}^{k}}(x+r).
\label{BIN2}%
\end{equation}
In order to prove this claim it is enough to show \eqref{BIN2} for
$x=-k,-k+1,...,0.$ Since for any natural \( m \),
\[
\frac{1}{l!}\prod_{r=0}^{l-1}(r-m)=
\begin{cases}
  (-1)^{l}\binom{m}{l},&\,l\leq m\\
  0,&\,l>m,
\end{cases}
\]
 we have to show that for \( m\leq k \),
\[
\sum_{l=1}^{m}\binom{m}{l}(-1)^{l}\left[  \sum_{j=1}^{l}\binom{l}{j}%
(-1)^{j}j^{k+1}\right]  =m^{k+1}.
\]
For $k=0$ the above equality is obvious.
Assume that we have proved the claim for $k\leq n$. Then it follows that
\begin{multline*}
\sum_{l=1}^{m}\binom{m}{l}(-1)^{l}\left[  \sum_{j=1}^{l}\binom{l}{j}%
(-1)^{j}j^{n+1}\right]  \\
  =m\sum_{l=1}^{m-1}\binom{m-1}{l}(-1)^{l}\left[
\sum_{j=1}^{l}\binom{l}{j}(-1)^{j}(j+1)^{n}\right]
\\
=m\sum_{s=0}^{n}\binom{n}{s}\sum_{l=1}^{m-1}\binom{m-1}{l}(-1)^{l}\left[
\sum_{j=1}^{l}\binom{l}{j}(-1)^{j} j^{s}\right]
\\
 =m\sum_{s=0}^{n}\binom{n}{s}(m-1)^{s}=m^{n+1}.
\end{multline*}
\end{proof}

\subsection*{Proof of Theorem~\ref{SG_AV}}

From the Taylor formula for semi-groups it
follows that%
\[
P_{s}f_{u}=\sum_{k=0}^{r}\frac{s^{k}}{k!}A^{k}f_{u}+\frac{1}{r!}\int_{0}%
^{s}(s-\tau)^{r}P_{\tau}A^{r+1}f_{u}d\tau.
\]
Due to (\ref{HG}), for $0\leq\tau\leq s<R_{u}$ we have%
\[
\left\Vert P_{\tau}A^{r+1}f_{u}\right\Vert \leq\sup_{0\leq\tau\leq
s}\left\Vert P_{\tau}\right\Vert \,\left\Vert A^{r+1}f_{u}\right\Vert \leq
\sup_{0\leq\tau\leq s}\left\Vert P_{\tau}\right\Vert \,\left(  \frac{1}{R_{u}%
}+\varepsilon\right)  ^{r+1}(r+1)!
\]
for any $\varepsilon>0.$ It thus follows that
\[
\left\Vert P_{s}f_{u}-\sum_{k=0}^{r}\frac{s^{k}}{k!}A^{k}f_{u}\right\Vert
\leq\left(  \frac{1}{R_{u}}+\varepsilon\right)  ^{r+1}s^{r+1}\sup_{0\leq\tau\leq
s}\left\Vert P_{\tau}\right\Vert,
\]
which converges to zero when $r\rightarrow\infty,$ if $\left\vert s\right\vert
<R_{u}/(1+\varepsilon R_{u}).$ Since $\varepsilon>0$ is arbitrary, the first statement
is proved.
\par
The commutation property $A^{k}P_{t}f_{u}=P_{t}A^{k}f_{u}$ and
the boundedness of $P_{t}$ for $t\geq0$ imply that for $\left\vert
s\right\vert <R_{u}$,%
\[
\sum_{k=0}^{\infty}\frac{\left\vert s\right\vert ^{k}}{k!}\left\Vert
A^{k}P_{t}f_{u}\right\Vert =\sum_{k=0}^{\infty}\frac{\left\vert s\right\vert
^{k}}{k!}\left\Vert P_{t}A^{k}f_{u}\right\Vert \leq\left\Vert P_{t}\right\Vert
\sum_{k=0}^{\infty}\frac{\left\vert s\right\vert ^{k}}{k!}\left\Vert
A^{k}f_{u}\right\Vert <\infty.
\]
Since $P_{t}f_{u}\in\mathfrak{D}(A^{k})$ for all $k\geq0$, (\ref{HG1})
follows.

\subsection*{Proof of Theorem~\ref{affthm}}
For $r\geq0$ define $A^{r}f_{u}$ $=:$ $g_{r}f_{u}$ with $f_{u}(x)=$
$\exp\left[  \mathfrak{i}u^{\top}x\right]  ,$ and write
\[
A^{r+1}f_{u}=A\left(  g_{r}\exp(\mathfrak{i}u^{\top}x)\right)  =\sum
_{|\alpha|\geq1}\mathfrak{a}_{\alpha}(x)\partial_{x^{\alpha}}\left(  g_{r}%
\exp(\mathfrak{i}u^{\top}x)\right)  .
\]
Leibniz formula implies
\begin{align*}
A^{r+1}f_{u}  & =\sum_{|\alpha|\geq1}\mathfrak{a}_{\alpha}(x)\sum_{\beta
\leq\alpha}\frac{\alpha!}{\beta!(\alpha-\beta)!}\partial_{x^{\beta}}%
g_{r}\partial_{x^{\alpha-\beta}}\exp(\mathfrak{i}u^{\top}x)\\
& =\left(  \sum_{|\alpha|\geq1}\mathfrak{a}_{\alpha}(x)\sum_{\beta\leq\alpha
}\frac{\alpha!}{\beta!(\alpha-\beta)!}(\mathfrak{i}u)^{\alpha-\beta}%
\partial_{x^{\beta}}g_{r}\right)  \exp(\mathfrak{i}u^{\top}x).
\end{align*}
Hence, the following recurrent formula holds
\begin{equation}
g_{r+1}=\sum_{|\alpha|\geq1}\mathfrak{a}_{\alpha}(x)\sum_{\beta\leq\alpha
}\frac{\alpha!}{\beta!(\alpha-\beta)!}\left(  \mathfrak{i}u\right)
^{\alpha-\beta}\partial_{x^{\beta}}g_{r}.\label{gp1}%
\end{equation}
Similar formulas for derivatives of $g_{r+1}$ can be obtained:
\[
\partial_{x^{\rho}}g_{r+1}=\sum_{|\alpha|\geq1}\sum_{\beta\leq\alpha}%
\frac{\alpha!}{\beta!(\alpha-\beta)!}\left(  \mathfrak{i}u\right)
^{\alpha-\beta}\sum_{\eta\leq\rho}\frac{\rho!}{\eta!(\rho-\eta)!}%
\partial_{x^{\eta}}\mathfrak{a}_{\alpha}\ \partial_{x^{\rho-\eta+\beta}}g_{r}.
\]
Since the underlying process is affine, all derivatives of $\mathfrak{a}$ of
order higher than one are zero and thus, by induction, $g_{r}$ is polynomial
in $x$ of degree at most equal $r.$ We so get for $\left\vert \rho\right\vert \leq r+1,$
\[
\partial_{x^{\rho}}g_{r+1}=\sum_{\substack{\eta\leq\rho,\,|\eta|\leq1}}\frac{\rho!}{\eta!(\rho-\eta)!}\sum_{|\alpha|\geq1}\sum_{\beta\leq\alpha}\frac{\alpha!}{\beta!(\alpha-\beta)!}\left(  \mathfrak{i}u\right)
^{\alpha-\beta}\partial_{x^{\eta}}\mathfrak{a}_{\alpha}\ \partial
_{x^{\rho-\eta+\beta}}g_{r}%
\]
By defining
\[
\Gamma_{r}:=\max_{|\beta|\leq r}\left\vert \partial_{x^{\beta}}g_{r}%
\right\vert,
\]
we obtain the following estimate for $x\in\mathfrak{X},$
\[
\left\vert \partial_{x^{\rho}}g_{r+1}\right\vert \leq\Gamma_{r}(1+\|x\|)\sum
_{\substack{\eta\leq\rho,\\|\eta|\leq1}}\frac{\rho!}{\eta!(\rho-\eta)!}%
\sum_{|\alpha|\geq1}\sum_{\beta\leq\alpha}\frac{\alpha!}{\beta!(\alpha-\beta)!}\left\vert u\right\vert ^{\alpha
-\beta}\mathcal{D}_{|\alpha|}^{\mathfrak{a}}%
\]
with $\left\vert u\right\vert :=[\left\vert u_{1}\right\vert ,...,\left\vert
u_{n}\right\vert ].$ Hence, by the simple relations
\begin{align*}
\sum_{\{\eta:\left\vert \eta\right\vert \leq1\}}\frac{\rho!}{\eta!(\rho-\eta)!} &
=1+|\rho|,\qquad\sum_{\beta\leq\alpha}\frac{\alpha!}{\beta!(\alpha-\beta
)!}\left\vert u\right\vert ^{\alpha-\beta}\leq(1+\|u\|)^{|\alpha|},\\
\sum_{\alpha,|\alpha|=k}1 &  =\frac{(k+n-1)!}{k!(n-1)!}\leq2^{n+k},
\end{align*}
with $\Vert u\Vert=\max_{i=1,\ldots,n}|u_{i}|,$ we have
\begin{eqnarray}
\Gamma_{r+1}
&\leq& 2^{n}\Gamma_{r}(r+2)(1+\|x\|)\sum_{k\geq1} 2^k(1+\|u\|)^{k}\mathcal{D}_{k}^{\mathfrak{a}}\label{series}
\\
\nonumber
&=&2^{n}\Gamma_{r}(r+2)(1+\|x\|)\theta(\|u\|),
\end{eqnarray}
where the series in (\ref{series}) is convergent due to assumption
(\ref{momres}).  As a consequence, (\ref{coefest}) holds. $\Box $

\subsubsection*{Proof of Proposition~\ref{affprop}}
From (\ref{gp1}) we have with $\mathfrak{a}_{0}:=0,$
\begin{align*}
g_{r+1}
&  =\sum_{\alpha,\beta,\gamma\geq0}\mathfrak{a}_{\alpha}\frac{\alpha!}%
{\beta!(\alpha-\beta)!}(\mathfrak{i}u)^{\alpha-\beta}g_{r,\gamma}\frac
{\gamma!}{(\gamma-\beta)!}x^{\gamma-\beta}1_{\left\vert \gamma\right\vert \leq
r}1_{\beta\leq\alpha}1_{\beta\leq\gamma}\\
&  =\sum_{\beta,\gamma\geq0}g_{r,\gamma+\beta}\binom{\gamma+\beta}{\beta}
x^{\gamma}1_{\left\vert \gamma+\beta\right\vert \leq r}\mathfrak{b}_{\beta
}\\
&  =\sum_{\left\vert \gamma\right\vert \leq r}x^{\gamma}\sum_{\left\vert
\beta\right\vert \leq r-\left\vert \gamma\right\vert }g_{r,\gamma+\beta}%
\binom{\gamma+\beta}{\beta}\mathfrak{b}_{\beta}^{0}\\
&  +\sum_{\,\left\vert \gamma\right\vert \leq r}\sum_{\,\left\vert
\beta\right\vert \leq r-\left\vert \gamma\right\vert }g_{r,\gamma+\beta}%
\binom{\gamma+\beta}{\beta}\sum_{\kappa,\,\left\vert \kappa\right\vert
=1}\mathfrak{b}_{\beta,\kappa}^{1}x^{\gamma+\kappa}\\
&  =\sum_{\,\left\vert \gamma\right\vert \leq r+1}x^{\gamma}\sum_{\left\vert
\beta\right\vert \leq r-\left\vert \gamma\right\vert }g_{r,\gamma+\beta}%
\binom{\gamma+\beta}{\beta}\mathfrak{b}_{\beta}^{0}\\
&  +\sum_{\,\left\vert \gamma\right\vert \leq r+1}x^{\gamma}\sum_{\left\vert
\kappa\right\vert =1,\ \kappa\leq\gamma}\sum_{\beta,\,\left\vert
\beta\right\vert \leq r+1-\left\vert \gamma\right\vert }g_{r,\gamma
-\kappa+\beta}\binom{\gamma-\kappa+\beta}{\beta}\mathfrak{b}_{\beta,\kappa
}^{1},%
\end{align*}
where empty sums are to be interpret as zero.
The second recursion follows from $(r+1)h_{r+1}=\eta^{-1}%
\widetilde{h}_{r+1}+rh_{r}$ with $A(h_{r}f_{u})=\widetilde{h}_{r+1}f_{u}$ and
$\widetilde{h}_{r+1}$ computed via (\ref{grecu}) with $g_{r}$ replaced
by $h_{r}.$ $\Box$

\end{document}